\documentclass[10pt]{article}

\usepackage{amsfonts,amsmath,graphicx}
\usepackage[final]{epsfig}
\usepackage{graphics}
\usepackage{amsmath}
\usepackage{amsfonts}
\usepackage{latexsym}
\usepackage{amssymb}
\usepackage{graphicx}
\usepackage{epstopdf}
\newtheorem{lemma}{Lemma}[section]

\newtheorem{theorem}{Theorem}

\newcommand{\proofend}{$\Box$\bigskip}

\newcommand{\Q}{{\mathbb Q}}

\newcommand{\RP}{{\mathbb{RP}}}

\newcommand{\1}{{\bf 1}}
\newcommand{\2}{{\bf 2}}
\newcommand{\3}{{\bf 3}}
\newcommand{\4}{{\bf 4}}
\newcommand{\5}{{\bf 5}}
\newcommand{\6}{{\bf 6}}

\def\proof{\paragraph{Proof.}}

\begin{document}
\date{\today}

\title {Periodic  trajectories in the regular pentagon, II}
\author{Dmitry Fuchs,\thanks{
Department of Mathematics,
University of California, Davis, CA 95616, USA;
e-mail: \tt{fuchs@math.ucdavis.edu}
}
\ and Serge Tabachnikov\thanks{
Department of Mathematics,
Pennsylvania State University, University Park, PA 16802, USA;
e-mail: \tt{tabachni@math.psu.edu}
}
}

\maketitle

\section{Introduction and formulation of results} \label{intro}

In our recent paper \cite{DFT}, we studied periodic billiard trajectories in a regular pentagon and in the isosceles triangle with  with the angles $(\pi/5,\pi/5,3\pi/5)$. We provided a full computation  of the lengths of these trajectories, both geometric and combinatoric, and formulated some conjectures concerning symbolic periodic trajectories. The main goal of this article is to prove two of these conjectures.

Technically, the study of billiard trajectories in a regular pentagon is essentially equivalent  to the study of geodesics in the ``double pentagon," a translation surface obtained from two centrally symmetric copies of a regular pentagon by pairwise pasting the parallel sides. The result is a surface of genus 2 that has a flat structure inherited from the plane and a conical singularity. See \cite{HS,MT,Sm,Tab,Vo,Zo} for surveys of flat surfaces and rational polygonal billiards.  

\begin{figure}[hbtp]
\centering
\includegraphics[width=3in]{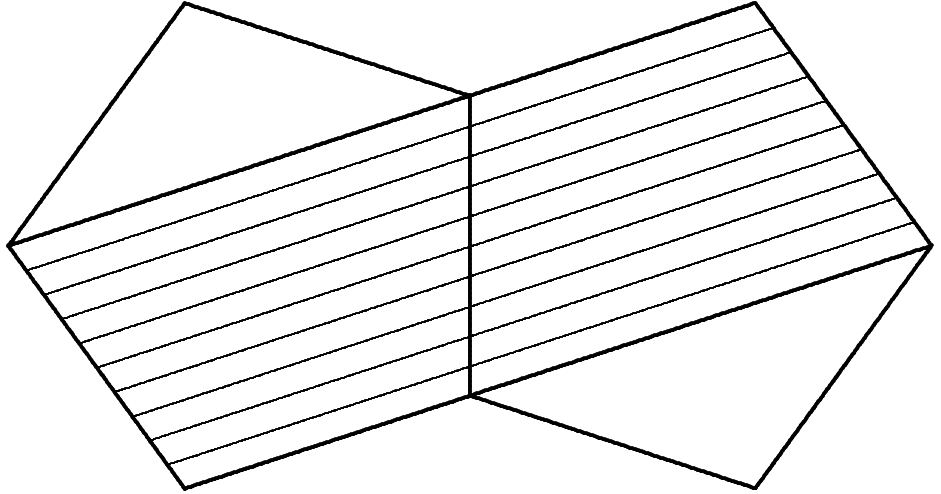}
\caption{Two strips of parallel periodic geodesics covering the double pentagon}
\label{strips}
\end{figure}

Let us describe the relevant results from \cite{DFT}. First of all, a periodic linear trajectory is always included into a parallel family of such trajectories, and when we talk about the period, length, symbolic orbit, etc., we always mean these parallel families. See Figure \ref{strips}. 

Second, the double pentagon has  an involution, the central symmetry that exchanges the two copies of the regular pentagon. This involution interchanges the linear trajectories that have the opposite directions. For this reason, we identify the opposite directions, so the set of directions is the real projective line $\RP^1$. We identify this projective line with the circle at infinity of the hyperbolic plane in the Poincar\'e disc model. 

It is clear from Figure \ref{strips} that the directions of the sides of the pentagons are periodic: every linear trajectory in this direction is closed. In fact, the so-called Veech dichotomy applies to the double pentagon: if there exists a periodic trajectory in some direction then all parallel trajectories are also periodic (and they form two strips, longer -- shaded in Figure \ref{strips}, and shorter -- left unshaded, covering the double pentagon). These periodic trajectories are constructed as follows.

\begin{figure}[hbtp]
\centering
\includegraphics[width=3.5in]{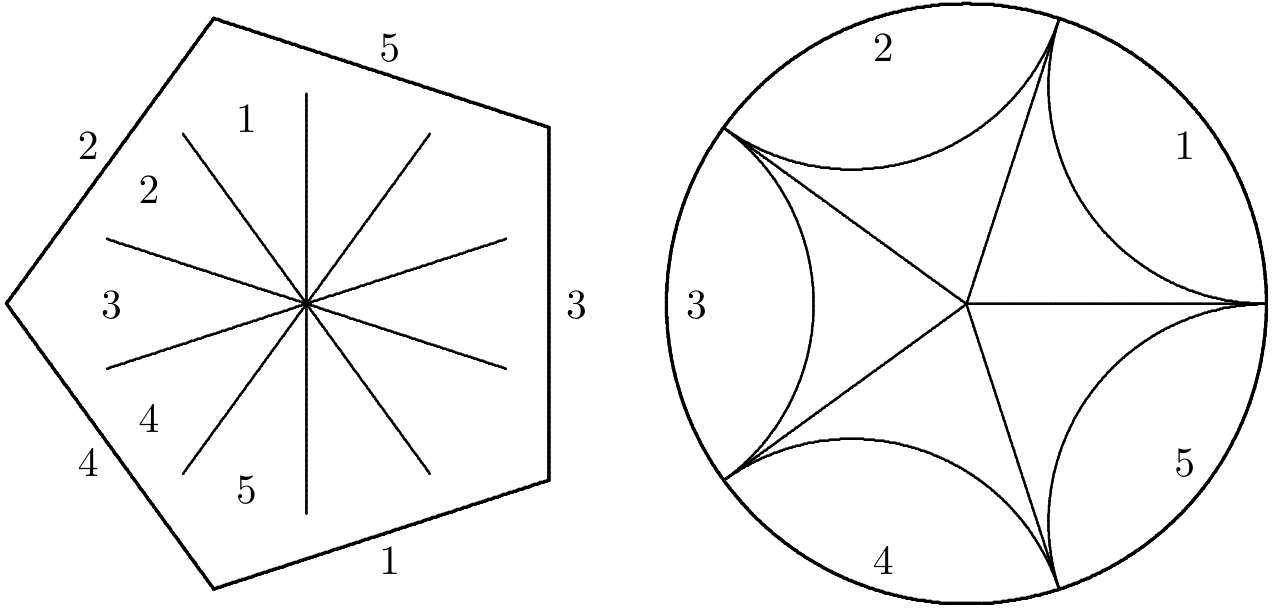}
\caption{Sectors of directions and the ideal pentagon}
\label{directions}
\end{figure}

Consider Figure \ref{directions}. On the left, one has the regular pentagon and the five periodic directions parallel to its sides. These directions form five cones, labeled 1 through 5. The same periodic directions on $\RP^1$ form a regular ideal pentagon in the hyperbolic plane; it is shown on the right. Due to rotational symmetry, we are interested in the periodic directions in the third sector (bounded, on the right part of the figure, by the vertical arc).

\begin{figure}[hbtp]
\centering
\includegraphics[width=4in]{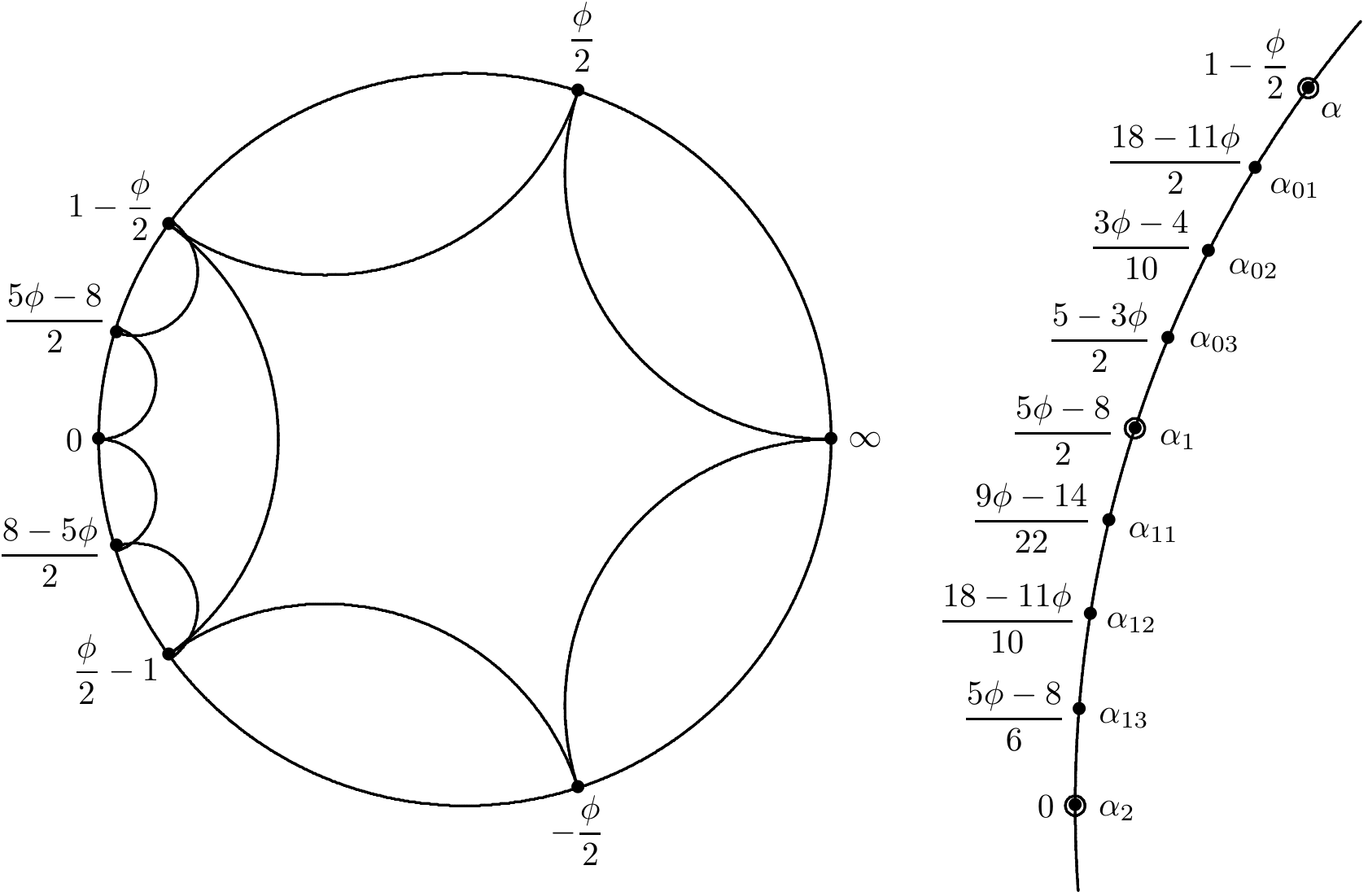}
\caption{Generating periodic directions: the pentagons of 0th and 1st generations}
\label{tiling}
\end{figure}

Let us refer to the  ideal pentagon in Figure \ref{directions}, right, as the pentagon of 0th generation. This pentagon 
tiles the hyperbolic plane by reflections in its sides; the tiling is made of the pentagons of the first, second, etc., generations. The vertices of this tiling, and only them, are the periodic directions in the double pentagon. We choose an affine coordinate on $\RP^1$ so that the vertices of the ideal pentagon have the coordinates 
$\displaystyle{1-\frac{\phi}2,\frac{\phi}2,\infty,-\frac{\phi}2,\frac{\phi}2-1}$ 
where $\displaystyle{\phi=\frac{1+\sqrt{5}}{2}}$ is the Golden Ratio. Then the periodic directions form the set $\Q[\phi]\cup\{\infty\}$. See Figure \ref{tiling} left.

The periodic directions on the  arc  $\left(1-\frac{\phi}{2},\frac{\phi}{2}-1\right)$ which are  the vertices of the pentagons of $k$th generation are denoted by $\alpha$ with $k$ indices as shown in Figure \ref{tiling}, right. For example, the pentagon bounded by the arc $(\alpha,\alpha_1)$ has vertices $\alpha,\alpha_{01},\alpha_{02},\alpha_{03},\alpha_1$, and the pentagon bounded by the arc $(\alpha_{011},\alpha_{012})$ has vertices $\alpha_{011},\alpha_{0111},\alpha_{0112},$ $\alpha_{0113},\alpha_{012}$  (note that the last index is never zero).

We label the pairs of parallel sides of the two regular pentagons which form the double pentagon by the symbols 1,\ 2,\ 3,\ 4,\ 5, as indicated in Figure \ref{directions}, left. Thus we talk about symbolic  trajectories, periodic words in these five symbols, consisting of the symbols of the sides that are consecutively crossed by a periodic linear trajectory. The combinatorial period of a trajectory is the period of this word. 

As we mentioned earlier, to every periodic direction, there correspond two strips of parallel periodic trajectories that cover the surface; we denote their periods by a pair consisting of a script and capital letter, such as $(a,A)$, with $a\le A$. For example, the symbolic trajectories on Figure \ref{strips} are $25$ and $43$, and their periods are $(2,2)$.

The following theorem was proved in \cite{DFT}.

\begin{figure}[hbtp]
\centering
\includegraphics[width=2in]{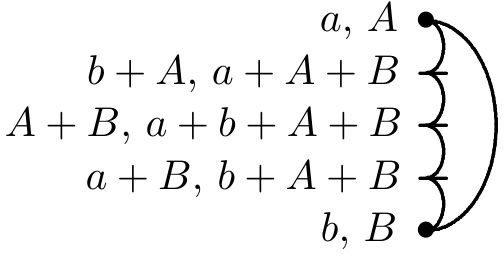}
\caption{The recurrence for periods}
\label{periods}
\end{figure}

\begin{theorem} \label{pers}
Let $(a,A)$ and $(b,B)$ be the pairs of periods for two points joined by a side of an ideal pentagon of some generation. Then, for the three additional vertices of the pentagon of the next generations, the periods are as shown in Figure \ref{periods}.
\end{theorem}

A different, yet equivalent, description of periods is given in the next theorem, also from \cite{DFT}.

\begin{figure}[hbtp]
\centering
\includegraphics[width=3in]{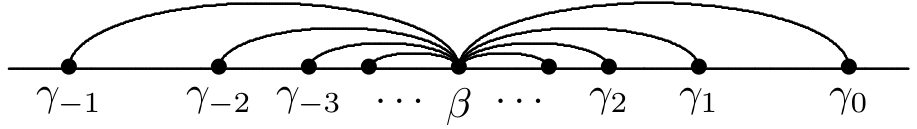}
\caption{Periodic directions connected with a given periodic direction}
\label{arcs}
\end{figure}

\begin{theorem} \label{persprime}
Let $\beta$ be a periodic direction, and let $\dots,\gamma_{-2},\gamma_{-1},$ $\gamma_0,\gamma_1,\gamma_2,\dots$ be all points connected by arcs (sides of ideal pentagons) with $\beta$ and ordered as shown in Figure \ref{arcs}. Let $(b,B)$ be the periods corresponding to the direction $\beta$, and let $(c_i,C_i)$ be the periods corresponding to $\gamma_i$, multiplied by $-1$, if $i<0$. Then 
$$\dots,(c_{-2},C_{-2}),(c_{-1},C_{-1}),(c_0,C_0),(c_1,C_1),(c_2,C_2),\dots$$
is an arithmetic progression with the difference $(B,b+B)$.
\end{theorem}

We observed, in computer experiments, that the recurrence relations of Theorem \ref{pers} and \ref{persprime} also held for the symbolic trajectories, and we conjectured the following results in \cite{DFT} which we formulate as the next two theorems.

\begin{theorem} \label{first} 
{\bf (Conjecture 8 of \cite{DFT}.)} Let $a,A$ and $b,B$ be two pairs of cyclic symbolic orbits, corresponding to two periodic directions  joined by a side of an ideal pentagon of some generation. Then one can cut the cyclic words $a,A,b,B$ into linear ones, concatenate them, and close the words up to cyclic words, so that 
the cyclic symbolic orbits for  the three additional vertices of the pentagon of the next generations (listed in the direction from the first point to the second)  are:
$$\begin{array} {ll} \ bA,&BaA;\\ A B, &bBaA;\\ aB,&AbB.\end{array}$$
\end{theorem}

\begin{figure}[hbtp]
\centering
\includegraphics[width=5in]{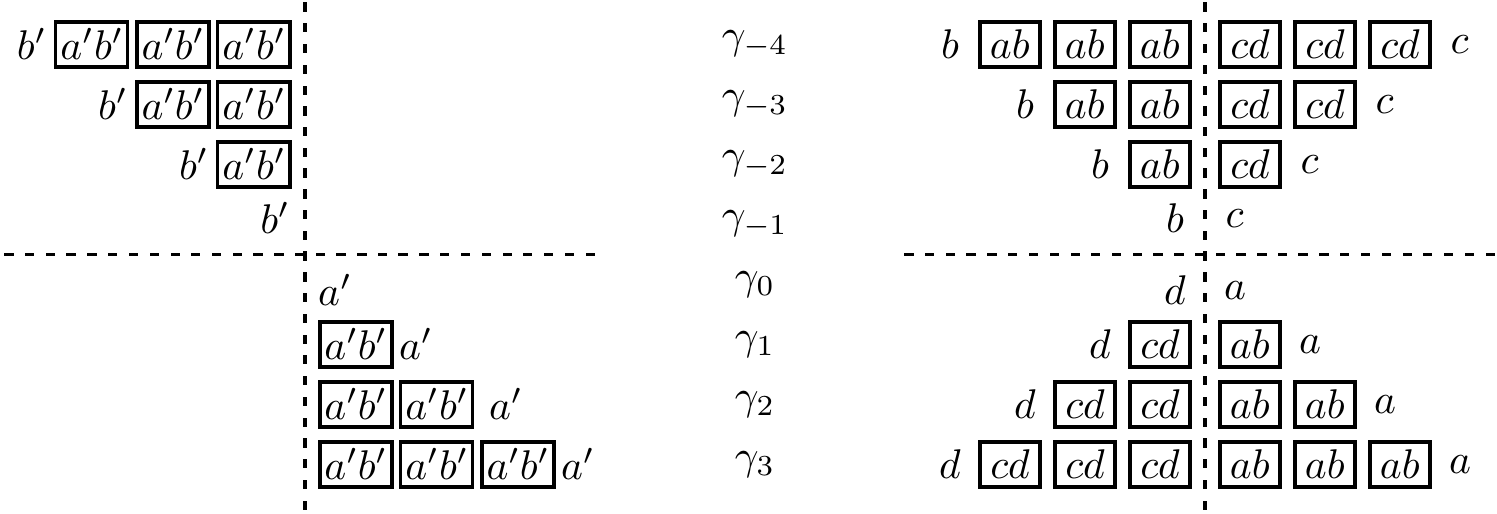}
\caption{Theorem \ref{second}}
\label{orbits}
\end{figure}

\begin{theorem} \label{second}
{\bf (Conjecture 9 of \cite{DFT}.)} Let $\beta$ and $\dots,\gamma_{-2},\gamma_{-1},\gamma_0,\gamma_1,\gamma_2,\dots$ denote the same as in Theorem \ref{persprime}. Then there exist a splitting of the short symbolic orbit corresponding to $\beta$ into two parts, $(c,d)$, and two splittings of the long symbolic orbit corresponding to $\beta$ into two parts: $(a,b)=(a',b')$ such that the short and long symbolic orbits corresponding to $\gamma_i$ are as shown in Figure \ref{orbits}. 
\end{theorem}
 
Notice that actually Conjecture 9 of \cite{DFT} contained one additional claim, namely that $a$ and $b'$ have the same beginning. This statement remains a conjecture.

We give proofs of Theorems \ref{first} and \ref{second} in the next two sections.  

\begin{figure}[hbtp]
\centering
\includegraphics[width=1.6in]{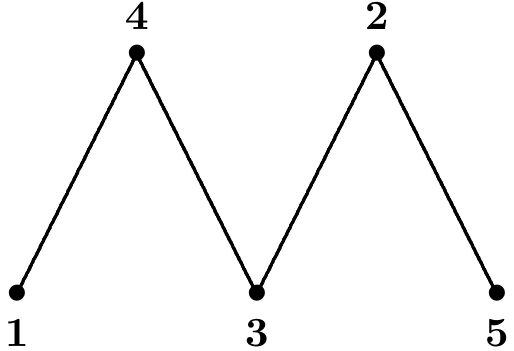}
\caption{Generating graph of symbolic orbits}
\label{graph}
\end{figure}

The main tool in the proofs is another result from \cite{DFT}, an algorithm producing symbolic periodic orbits of the next generation from the previous one. It is based on the graph shown in Figure \ref{graph}.  

Consider a symbolic orbit $W$ corresponding to $\alpha_{n_1n_2\dots n_k}$. For $j=1,2,3,4$, define the transformation $T_j$ as follows. Subtract $j$ mod 5 from each entry of $W$ to obtain a new periodic word in symbols 1,\ 2,\ 3,\ 4,\ 5. This word determines a  path on the graph in Figure \ref{graph}. When going along this path, we insert the symbols that we pass into the cyclic word; the resulting periodic word is $T_j(W)$. For example, $T_1(2343)=(14323234)$. 

One has the following result from \cite{DFT}.

\begin{theorem} \label{symb}
Let $W$ be the symbolic orbit, long or short, corresponding to the periodic direction  $\alpha_{n_1n_2\dots n_k}$. Then $T_j(W)$ is the symbolic orbit, long or short, respectively, corresponding to the periodic direction $\alpha_{j-1,3-n_1,3-n_2,\dots,3-n_{k-1},4-n_k}$. 
\end{theorem}

\section{Proof of Theorem \ref{first}}

In this section we prove Theorem \ref{first}. Let us summarize the set-up, as described in the previous section. 

Each symbolic orbit is a cyclic sequence of even length in the alphabet $\{1,2,3,4,5\}$  consisting of pairs $43,41,25,23$. To every vertex of an ideal pentagon of each generation, there correspond two symbolic orbits, the short one and the long one. The statement of the theorem is that the symbolic orbits corresponding to two consecutive vertices of each pentagon can be cut into linear sequences $a,A$ and $b,B$ (also consisting of pairs  $43,41,25,23$, and thus starting with 4 or 2 and ending with 1, or 3, or 5) in such a way that the symbolic orbits, corresponding to the three intermediate vertices of the pentagon of the next generation, are $bA,AaB;\ AB,AaBb;\ aB,BbA$.

We shall apply the transformations $T_j,\ j=1,2,3,4,$ to {\sl linear} (not cyclic) sequences. This makes it necessary to specify two additional things. First,  we place the path from the last digit to the first digit after the last digit. Second, after applying $T_j$, a cyclic permutation may be needed; we {\sl always} do it in the following way: if the first digit is 3, then we move it to the end of the sequence; if the first digit is not 3, but the last digit is 2  or 4, then we move this 2 or 4 to the beginning of the sequence.

The proof of the theorem is by induction on the generation of an ideal pentagon. 
\medskip

{\bf Base of induction}.
For the ``initial" pentagon (0th generation), the symbolic orbits are
$$a=41,A=23;\ \ bA,AaB;\ \ AB,AaBb;\ \ aB,BbA;\ \ b=25,B=43.$$
In particular, the statement of the theorem holds.
\medskip

{\bf Induction setting}.
We assume that, for some pentagon of some generation, the statement of the theorem holds. The symbolic orbits of the two extreme vertices are of the form $a=x_1\ast x_2,\ A=X_1\ast X_2,\ b=y_1\ast y_2,\ B=Y_1\ast Y_2$ where $x_1,\dots,Y_2 \in \{1,2,3,4,5\}$. One has nine cases to consider:
$$
\begin{matrix}
({\1})&a=4\ast1&A=2\ast3&b=2\ast5&B=4\ast3\\
 ({\2})&a=2\ast3&A=4\ast3&b=4\ast1&B=2\ast3\\ 
 ({\3})&a=2\ast3&A=2\ast3&b=2\ast3&B=2\ast3\\
  ({\4})&a=4\ast3&A=4\ast3&b=4\ast3&B=4\ast3\\
   ({\5})&a=2\ast5&A=4\ast3&b=4\ast3&B=2\ast3\\
    ({\6})&a=4\ast3&A=2\ast3&b=2\ast3&B=4\ast3\\
\end{matrix}
$$
and also $({\2})^\ast,({\5})^\ast,({\6})^\ast$ where $\ast$ means the involution $(a,A)\leftrightarrow(b,B)$.

We will prove that after every transformation $T_j$, we  get a pentagon satisfying the statement of the theorem, with $a,A,b,B$ satisfying the same conditions; moreover, the beginnings and the ends of the new $a,A,b,B$ are determined by the  sequence $(\1)-(\6)^\ast$. Specifically, the transformations $T_j$ act in the set $(\1),\dots,(\6)^\ast$ in the following way:
$$
\begin{matrix}
T_1(\1)=(\2)&T_2(\1)=(\3)&T_3(\1)=(\4)&T_4(\1)=(\5)\\
T_1(\2)=(\6)&T_2(\2)=(\3)&T_3(\2)=(\4)&T_4(\2)=(\5)^\ast\\
T_1(\3)=(\4)&T_2(\3)=(\3)&T_3(\3)=(\4)&T_4(\3)=(\4)\\
T_1(\4)=(\3)&T_2(\4)=(\3)&T_3(\4)=(\4)&T_4(\4)=(\3)\\
T_1(\5)=(\2)^\ast&T_2(\5)=(\3)&T_3(\5)=(\4)&T_4(\5)=(\6)\\
T_1(\6)=(\6)^\ast&T_2(\6)=(\3)&T_3(\6)=(\4)&T_4(\6)=(\6)^\ast\\
\end{matrix}
$$
and, obviously, $T_j({\bf k})^\ast=(T_j({\bf k}))^\ast$. 

Since the initial pentagon corresponds to the case (\1) (with empty $\ast$), this will prove our statement.
\medskip

{\bf Induction step}. 
The rest of the proof is a tedious direct computation. To make it slightly less tedious, we introduce some  notations and make a preliminary remark.
\smallskip

{\bf Preparation}.
When we apply $T_j$ to one of the sequences $a,A,bA,AaB,\\
\dots,B$, we first apply it to each symbol of $a,A,b,B$, and then insert the digits of the connecting path in parentheses. For example, $T_2(4\ast1\, 4\ast 3)=2\ast4(3)2\ast1(43)$. 

Next, if we need to apply a cyclic permutation, say, to move 3 to the end of the word, we draw a right arrow above the digit 3 in the beginning of the word and add  3 in angular brackets in the end (which  shows that the first 3 is not there, but the last 3 is there). The same with movable 2 and 4. For example, $T_1(2\ast3)=\langle4\rangle1\ast2(3\overleftarrow4)$ (the latter means $41\ast23$). 

And the last thing. In many cases, the first or the last digit of a sequence provides an information about the second digit or the digit before the last after applying some $T_j$; in some cases this is important for us. For example, consider $T_1(4\ast z)$ (what $z$ is does not matter to us). The sequence $4\ast z$ starts with $43$ or $41$. After applying $T_1$, we get $32$ or $35$; after inserting the path between 3 and 5, the second word becomes $325$. We see that $T_1(4\ast z)$ starts not just with $3$, but rather with $32$. We reflect this in our notations: we write, for example, $T_1(4\ast1)=\overrightarrow32\ast5(2)\langle3\rangle$ (not just $T_1(4\ast1)=\overrightarrow3\ast5(2)\langle3\rangle$. There are several other cases similar to this; we prefer to list them in the next lemma.

\begin{lemma} \label{lm} One has: 
\begin{equation*}
\begin{split}
T_1(4\ast\dots)=32\ast\dots,\ T_4(2\ast\dots)=34\ast\dots,\\
T_1(\dots\ast5)=\dots\ast14,\  T_4(\dots\ast3)=\dots\ast34,\ T_2(\dots\ast1)=\dots\ast34, \\
T_3(\dots\ast5)=\dots\ast32,\  T_1(\dots\ast3)=\dots\ast32,\ T_4(\dots\ast1)=\dots\ast52.
\end{split}
\end{equation*}
\end{lemma}
\medskip

{$T_j({\bf k})$ {\bf for} ${\bf k}=\3,\4$.}
We need to apply 4 transformations in 6 cases. Since for ${\bf k}=\3,\4$, when all $a,A,b,B$ are the same, this requires much less work, we will start with these cases. In our notations, 
$$
\begin{matrix}
T_1(2\ast3)=\langle4\rangle1\ast32(3\overleftarrow4);&T_1(4\ast3)=\overrightarrow32\ast32\langle3\rangle;\\
 T_2(2\ast3)=\langle2\rangle5\ast2(43\overleftarrow2);&T_2(4\ast3)=2\ast1(43);\hfil\\
  T_3(2\ast3)=4\ast5(23);&T_3(4\ast3)=\langle4\rangle1\ast5(23\overleftarrow4);\\
   T_4(2\ast3)=\overrightarrow34\ast34\langle3\rangle;&T_4(4\ast3)=\langle2\rangle5\ast4(3\overleftarrow2);\\ 
\end{matrix}
$$
which is compatible with the description of $T_j(\3),T_j(\4)$ given above. Also, the statement of the theorem holds: the symbolic orbits for the intermediate vertices are simply the 2-, 3-, or 4-fold repetitions of the symbolic orbits for the extreme vertices, and applying $T_j$ will lead to similar repetitions. For example,
$$
T_1(2\ast3\, 2\ast3)=\langle4\rangle1\ast32(34)1\ast32(3\overleftarrow4),
$$
which is a twofold repetition of $41\ast323$.
\medskip

{ $T_j({\bf k})$ {\bf for} ${\bf k}=\1,\2,\5,\6$.}
Here we cannot avoid applying all the $T_j$ to all the symbolic orbits in all the four cases. It is done below, but necessarily occupies $4\times5\times4=80$ lines. To make  the validity of the statement of the theorem more visible, we separate the symbolic orbits for the intermediate vertices into pieces, corresponding to $a,A,b,B$, by the symbol $\|$.
\medskip

\begin{tabular}{lll}
&$\overrightarrow32\ast5(2)\langle3\rangle$&$\langle4\rangle1\ast32(3\overleftarrow4)$\\
\hskip-.3in$T_1(\1)=(\2)$&$\langle4\rangle1\ast1\|41\ast32(3\overleftarrow4)$&$\langle4\rangle1\ast323\|2\ast5(2)3\|2\ast32(3\overleftarrow4)$\\
&$\langle4\rangle1\ast323\|2\ast32(3\overleftarrow4)$&$\langle4\rangle1\ast323\|2\ast5(2)3\|2\ast32(3\|4)1\ast1\overleftarrow4$\\
&$\overrightarrow32\ast5(2)3\|2\ast32\langle3\rangle$&$\overrightarrow32\ast32(3\|4)1\ast1\|41\ast32\langle3\rangle$\\
&$\langle4\rangle1\ast1\overleftarrow4$&$\overrightarrow32\ast32\langle3\rangle$\\
\end{tabular}
\smallskip
\hrule
\smallskip
\begin{tabular}{lll}
&\hskip-.08in$2\ast34(3)$&\hskip-.1in$\langle2\rangle5\ast1(43\overleftarrow2)$\\ \hskip-.3in$T_2(\1)=(\3)$&\hskip-.1in$\langle2\rangle5\ast3\|(2)5\ast1(43\overleftarrow2)$&\hskip-.1in$\langle2\rangle5\ast1(43)\|2\ast34(3)\|2\ast1(43\overleftarrow2)$\\
&\hskip-.1in$\langle2\rangle5\ast1(43)\|2\ast1(43\overleftarrow2)$&\hskip-.1in$\langle2\rangle5\ast1(43)\|2\ast34(3)\|2\ast1(43\|2)5\ast3(\overleftarrow2)$\\
&\hskip-.08in$2\ast34(3)\|2\ast1(43)$&\hskip-.08in$2\ast1(43\|2)5\ast3\|(2)5\ast1(43)$\\
&\hskip-.1in$\langle2\rangle5\ast3(\overleftarrow2)$&\hskip-.08in$2\ast1(43)$\\
\end{tabular}
\smallskip
\hrule
\smallskip
\begin{tabular}{lll}
&$\langle4\rangle1\ast3(\overleftarrow4)$&$4\ast5(23)$\\
\hskip-.3in$T_3(\1)=(\4)$&$4\ast32(3)\|4\ast5(23)$&$4\ast5(23\|4)1\ast3\|(4)1\ast5(23)$\\
&$4\ast5(23\|4)1\ast5(23)$&$4\ast5(23\|4)1\ast3\|(4)1\ast5(23)\|4\ast32(3)$\\
&$\langle4\rangle1\ast3\|(4)1\ast5(23\overleftarrow4)$&$\langle4\rangle1\ast5(23)\|4\ast32(3)\|4\ast5(23\overleftarrow4)$\\
&$4\ast32(3)$&$\langle4\rangle1\ast5(23\overleftarrow4)$\\
\end{tabular}
\smallskip
\hrule
\smallskip
\begin{tabular}{lll}
&$\langle2\rangle5\ast5\overleftarrow2$&\hskip.1in$\overrightarrow34\ast34\langle3\rangle$\\
\hskip-.3in$T_4(\1)=(\5)$&$\overrightarrow34\ast1\|4\ast34\langle3\rangle$&\hskip.1in$\overrightarrow34\ast34(3\|2)5\ast5\|25\ast4\langle3\rangle$\\
&$\langle2\rangle5\ast5\|25\ast4(3\overleftarrow2)$&\hskip.1in$\overrightarrow34\ast34(3\|2)5\ast5\|25\ast43\|4\ast1(4)\langle3\rangle$\\
&$\langle2\rangle5\ast5\|25\ast4(3\overrightarrow2)$&\hskip.1in$\langle2\rangle5\ast43\|4\ast1(4)3\|4\ast34(3\overrightarrow2)$\\
&$\overrightarrow34\ast1(4)\langle3\rangle$&\hskip.1in$\langle2\rangle5\ast4(3\overleftarrow2)$\\
\end{tabular}
\smallskip
\hrule
\smallskip
\begin{tabular}{lll}
&\hskip-.1in$\langle4\rangle1\ast32(3\overleftarrow4)$&\hskip-.1in$\overrightarrow32\ast32\langle3\rangle$\\
\hskip-.3in$T_1(\2)=(\6)$&\hskip-.12in$\overrightarrow32\ast53\|2\ast32\langle3\rangle$&\hskip-.1in$\overrightarrow32\ast32(3\|4)1\ast32(3\|4)1\ast32\langle3\rangle$\\
&\hskip-.1in$\overrightarrow32\ast32(3\|4)1\ast32\langle3\rangle$&\hskip-.1in$\overrightarrow32\ast32(3\|4)1\ast32(3\|4)1\ast323\|2\ast5\langle3\rangle$\\
&\hskip-.1in$\langle4\rangle1\ast32(3\|4)1\ast32(3\overleftarrow4)$&\hskip-.1in$\langle4\rangle1\ast323\|2\ast53\|2\ast32(3\overleftarrow4)$\\
&\hskip-.1in$\overrightarrow32\ast5\langle3\rangle$&\hskip-.1in$\langle4\rangle1\ast32(3\overleftarrow4)$\\
\end{tabular}
\smallskip
\hrule
\smallskip
\begin{tabular}{lll}
&\hskip-.12in$\langle2\rangle5\ast1(43\overleftarrow2)$&\hskip-.1in$2\ast1(43)$\\
\hskip-.3in$T_2(\2)=(\3)$&\hskip-.1in$2\ast34(3)\|2\ast1(43)$&\hskip-.1in$2\ast1(43\|2)5\ast1(43\|2)5\ast1(43)$\\
&\hskip-.1in$2\ast1(43\|2)5\ast1(43)$&\hskip-.1in$2\ast1(43\|2)5\ast1(43\|2)5\ast1(43)\|2\ast34(3)$\\
&\hskip-.12in$\langle2\rangle5\ast1(43\|2)5\ast1(43\overleftarrow2)$&\hskip-.12in$\langle2\rangle5\ast1(43)\|2\ast34(3)\|2\ast1(43\overleftarrow2)$\\
&\hskip-.1in$2\ast34(3)$&\hskip-.12in$\langle2\rangle5\ast1(43\overleftarrow2)$\\
\end{tabular}
\smallskip
\hrule
\smallskip
\begin{tabular}{lll}
&\hskip-.08in$4\ast5(23)$&\hskip-.1in$\langle4\rangle1\ast5(23\overleftarrow4)$\\
\hskip-.3in$T_3(\2)=(\4)$&\hskip-.1in$\langle4\rangle1\ast3\|(4)1\ast5(23\overleftarrow4)$&\hskip-.1in$\langle4\rangle1\ast5(23)\|4\ast5(23)\|4\ast5(23\overleftarrow4)$\\
&\hskip-.1in$\langle4\rangle1\ast5(23)\|4\ast5(23\overleftarrow4)$&\hskip-.12in$\langle4\rangle1\ast5(23)\|4\ast5(23)\|4\ast5(23\|4)1\ast3(\overleftarrow4)$\\
&\hskip-.08in$4\ast5(23)\|4\ast5(23)$&\hskip-.08in$4\ast5(23\|4)1\ast3\|(4)1\ast5(23)$\\
&\hskip-.1in$\langle4\rangle1\ast3(\overleftarrow4)$&\hskip-.08in$4\ast5(23)$\\
\end{tabular}
\smallskip
\hrule
\smallskip
\begin{tabular}{lll}
&$\overrightarrow34\ast4\langle3\rangle$&$\langle2\rangle5\ast4(3\overleftarrow2)$\\
\hskip-.3in$T_4(\2)=(\5)^\ast$&$\langle2\rangle5\ast3\|(2)5\ast4(3\overleftarrow2)$&$\langle2\rangle5\ast43\|4\ast43\|4\ast4(3\overleftarrow2)$\\
&$\langle2\rangle5\ast43\|4\ast4(3\overleftarrow2)$&$\langle2\rangle5\ast43\|4\ast43\|4\ast4(3\|2)5\ast3(\overleftarrow2)$\\
&$\overrightarrow34\ast43\|4\ast4\langle3\rangle$&$\overrightarrow34\ast4(3\|2)5\ast3\|4\ast4\langle3\rangle$\\
&$\langle2\rangle5\ast3\overleftarrow2$&$\overrightarrow34\ast4\langle3\rangle$\\
\end{tabular}
\smallskip
\hrule
\smallskip
\begin{tabular}{lll}
&$\langle4\rangle1\ast1\overleftarrow4$&$\overrightarrow32\ast32\langle3\rangle$\\
\hskip-.3in$T_1(\5)=(\2)^\ast$&$\overrightarrow32\ast323\|2\ast32\langle3\rangle$&$\overrightarrow32\ast32(3\|4)1\ast1\|41\ast32\langle3\rangle$\\
&$\overrightarrow32\ast32(3\|4)1\ast32\langle3\rangle$&$\overrightarrow32\ast32(3\|4)1\ast1\|41\ast323\|2\ast32\langle3\rangle$\\
&$\langle4\rangle1\ast1\|41\ast32(3\overleftarrow4)$&$\langle4\rangle1\ast323\|2\ast323\|2\ast32(3\overleftarrow4)$\\
&$\overrightarrow32\ast32\langle3\rangle$&$\langle4\rangle1\ast32(3\overleftarrow4)$\\
\end{tabular}
\smallskip
\hrule
\smallskip
\begin{tabular}{lll}
&$\langle2\rangle5\ast3(\overleftarrow2)$&$2\ast1(43)$\\
\hskip-.3in$T_2(\5)=(\3)$&$2\ast1(43)\|2\ast1(43)$&$2\ast1(43\|2)5\ast3\|(2)5\ast1(43)$\\
&$2\ast1(43\|2)5\ast1(43)$&$2\ast1(43\|2)5\ast3\|(2)5\ast1(43)\|2\ast1(43)$\\
&$\langle2\rangle5\ast3\|(2)5\ast1(43\overleftarrow2)$&$\langle2\rangle5\ast1(43)\|2\ast1(43)\|2\ast1(43\overleftarrow2)$\\
&$2\ast1(43)$&$\langle2\rangle5\ast1(43\overleftarrow2)$\\
\end{tabular}
\smallskip
\hrule
\smallskip
\begin{tabular}{lll}
\hskip-.3in$T_3(\5)=(\4)$&\hskip-.08in$4\ast32(3)$&\hskip-.1in$\langle4\rangle1\ast5(23\overleftarrow4)$\\
&\hskip-.3in$\langle4\rangle1\ast5(23\|4)1\ast5(23\overleftarrow4)$&\hskip-.1in$\langle4\rangle1\ast5(23)\|4\ast32(3)\|4\ast5(23\overleftarrow4)$\\
&\hskip-.3in$\langle4\rangle1\ast5(23)\|4\ast5(23\overleftarrow4)$&\hskip-.1in$\langle4\rangle1\ast5(23)\|4\ast32(3)\|4\ast5(23\|4)1\ast5(23\overleftarrow4)$\\
&\hskip-.3in$4\ast32(3)\|4\ast5(23)$&\hskip-.1in$4\ast 5(23\|4)1\ast5(23\|4)1\ast5(23)$\\
&\hskip-.3in$\langle4\rangle1\ast5(23\overleftarrow4)$&\hskip-.1in$4\ast5(23)$\\
\end{tabular}
\smallskip
\hrule
\smallskip
\begin{tabular}{lll}
&\hskip-.1in$\overrightarrow34\ast1(4)\langle3\rangle$&\hskip-.1in$\langle2\rangle5\ast14(3\overleftarrow2)$\\
\hskip-.3in$T_4(\5)=(\6)$&\hskip-.1in$\langle2\rangle5\ast14(3\|25\ast14(3\overleftarrow2)$&\hskip-.1in$\langle2\rangle5\ast143\|4\ast1(4)3\|4\ast1(43\overleftarrow2)$\\
&\hskip-.12in$\langle2\rangle5\ast143\|4\ast1(43\overleftarrow2)$&\hskip-.14in$\langle2\rangle5\ast143\|4\ast1(4)3\|4\ast1(43\|2)5\ast14(3\overleftarrow2)$\\
&\hskip-.1in$\overrightarrow34\ast1(4)\|34\ast1(4)\langle3\rangle$&\hskip-.1in$\overrightarrow34\ast1(43\|2)5\ast14(3\|2)5\ast14\langle3\rangle$\\
&\hskip-.1in$\langle2\rangle5\ast14(3\overleftarrow2)$&\hskip-.1in$\overrightarrow34\ast1(4)\langle3\rangle$\\
\end{tabular}
\smallskip
\hrule
\smallskip
\begin{tabular}{lll}
&\hskip-.12in$\overrightarrow32\ast32\langle3\rangle$&\hskip-.12in$\langle4\rangle1\ast32(3\overleftarrow4)$\\ 
\hskip-.3in$T_1(\6)=(\6)^\ast$&\hskip-.12in$\langle4\rangle1\ast32(3\|4)1\ast32(3\overleftarrow4)$&\hskip-.12in$\langle4\rangle1\ast323\|2\ast323\|2\ast32(3\overleftarrow4)$\\
&\hskip-.14in$\langle4\rangle1\ast323\|2\ast32(3\overleftarrow4)$&\hskip-.14in$\langle4\rangle1\ast323\|2\ast323\|2\ast32(3\|4)1\ast32(3\overleftarrow4)$\\
&\hskip-.12in$\overrightarrow32\ast323\|2\ast32\langle3\rangle$&\hskip-.12in$\overrightarrow32\ast32(3\|4)1\ast32(3\|4)1\ast32\langle3\rangle$\\
&\hskip-.12in$\langle4\rangle1\ast32(3\overleftarrow4)$&\hskip-.12in$\overrightarrow32\ast32\langle3\rangle$\\
\end{tabular}
\smallskip
\hrule
\smallskip
\begin{tabular}{lll}
\hskip-.3in$T_2(\6)=(\3)$&\hskip-.1in$2\ast1(43)$&\hskip-.1in$\langle2\rangle5\ast1(43\overleftarrow2)$\\
&\hskip-.28in$\langle2\rangle5\ast1(43\|2)5\ast1(43\overleftarrow2)$&\hskip-.12in$\langle2\rangle5\ast1(43)\|2\ast1(43)\|2\ast1(43\overleftarrow2)$\\
&\hskip-.28in$\langle2\rangle5\ast1(43)\|2\ast1(43\overleftarrow2)$&\hskip-.18in$\langle2\rangle5\ast1(43)\|2\ast1(43)\|2\ast1(43\|2)5\ast1(43\overleftarrow2)$\\
&\hskip-.28in$2\ast1(43)\|2\ast1(43)$&\hskip-.08in$2\ast1(43\|2)5\ast1(43\|2)5\ast1(43)$\\
&\hskip-.28in$\langle2\rangle5\ast1(43\overleftarrow2)$&\hskip-.08in$2\ast1(43)$\\
\end{tabular}
\smallskip
\hrule
\smallskip
\begin{tabular}{lll}
&\hskip-.12in$\langle4\rangle1\ast5(23\overleftarrow4)$&\hskip-.1in$4\ast5(23)$\\
\hskip-.3in$T_3(\6)=(\4)$&\hskip-.08in$4\ast5(23)\|4\ast5(23)$&\hskip-.1in$4\ast5(23\|4)1\ast5(23\|4)1\ast5(23)$\\
&\hskip-.08in$4\ast5(23\|4)1\ast5(23)$&\hskip-.1in$4\ast5(23\|4)1\ast5(23\|4)1\ast5(23)\|4\ast5(23)$\\
&\hskip-.12in$\langle4\rangle1\ast5(23\|4)1\ast5(23\overleftarrow4)$&\hskip-.12in$\langle4\rangle1\ast5(23)\|4\ast5(23)\|4\ast5(23\overleftarrow4)$\\
&\hskip-.08in$4\ast5(23)$&\hskip-.12in$\langle4\rangle1\ast5(23\overleftarrow4)$\\
\end{tabular}
\smallskip
\hrule
\smallskip
\begin{tabular}{lll}
&\hskip-.12in$\langle2\rangle5\ast34(3\overleftarrow2)$&\hskip-.14in$\overrightarrow34\ast34\langle3\rangle$\\
\hskip-.3in$T_4(\6)=(\6)^\ast$&\hskip-.12in$\overrightarrow34\ast343\|4\ast34\langle3\rangle$&\hskip-.14in$\overrightarrow34\ast34(3\|2)5\ast34(3\|2)5\ast34\langle3\rangle$\\
&\hskip-.12in$\overrightarrow34\ast3(43\|2)5\ast34\langle3\rangle$&\hskip-.14in$\overrightarrow34\ast34(3\|2)5\ast34(3\|2)5\ast343\|4\ast34\langle3\rangle$\\
&\hskip-.12in$\langle2\rangle5\ast34(3\|2)5\ast34(3\overleftarrow2)$&\hskip-.12in$\langle2\rangle5\ast343\|4\ast343\|4\ast34(3\overleftarrow2)$\\
&\hskip-.1in$\overrightarrow34\ast34\langle3\rangle$&\hskip-.12in$\langle2\rangle5\ast34(3\overleftarrow2)$\\
\end{tabular}
\medskip

This finishes the proof.

\section{Proof of Theorem \ref{second}}

{\bf The statement}.  
 Let $\beta$ be a vertex of one of the pentagons, and let $\gamma_1,\gamma_2,\dots$ and $\gamma_{-1},\gamma_{-2},\dots$ be the two monotonic sequences of vertices of pentagons of further generations containing all vertices connected with $\beta$ by arcs of these pentagons, see Figure \ref{arcs}. The claim is that the symbolic orbits for $\gamma_n$ have the following form (we list first the short orbit and then the long one):
\begin{equation} \label{orbs}
\begin{split}
{\rm for}\ \gamma_{-1-n},n\ge0 \qquad\qquad A^nx, \qquad \overline A^nzva^n\\
{\rm for}\ \gamma_{1+n},n\ge0 \qquad\qquad\ yA^n,\qquad a^nuw\overline A^n
\end{split}
\end{equation}
where $a$ and $A$ are the symbolic orbits for $\beta$,  $\overline A$ is a cyclic permutation of $A$ (thus, $A$ and $\overline A$ present the same cyclic orbit, that is, the long orbit corresponding to $\beta$), $a=uv;\, A=xy;\, \overline A=zw$. In exceptional cases, some of $x,y,\dots,w$ will be empty.
\medskip

 {\bf Base of induction}. The proof  is by induction again, and we start with the pentagon of 0th generation. Its vertices are the points $\alpha_1,\alpha_2,\alpha_3$, and we need to check (\ref{orbs}) for these three values of $\beta$.
 
 \begin{lemma} \label{ind}
 Formulas (\ref{orbs}) hold with the following $x,y,\dots$:
 
 \begin{tabular} {lllllll}
 \hskip-.3in${\rm for}\ \alpha_1\colon$&\hskip-.15in$u=2523$&\hskip-.08in$v=\emptyset$&\hskip-.08in$x=4$&\hskip-.08in$y=4323$&$z=23$&$w=4143$\\
   \hskip-.3in&&&&\hskip-.1in$a=2523$&\hskip-.06in$A=414323$&\hskip-.04in$\overline A=234143$\\
   \hskip-.3in${\rm for}\ \alpha_2\colon$&\hskip-.15in$u=43$&\hskip-.1in$v=23$&\hskip-.08in$x=2523$&\hskip-.08in$y=4143$&$z=4143$&$w=2523$\\
   \hskip-.3in&&&&\hskip-.1in$a=4323$&\hskip-.06in$A=25234143$&\hskip-.04in$\overline A=41432523$\\
   \hskip-.3in${\rm for}\ \alpha_3\colon$&\hskip-.15in$u=4143$&\hskip-.08in$v=\emptyset$&\hskip-.08in$x=25$&\hskip-.08in$y=2343$&$z=43$&$w=2523$\\
   \hskip-.3in&&&&\hskip-.1in$a=4143$&\hskip-.06in$A=252343$&\hskip-.04in$\overline A=432523$\\ 
 \end{tabular}
 \end{lemma}
  
\proof
The symbolic orbits corresponding to the direction $\alpha$ are 41 and 25. Let $T_j$ be the same transformation as in Section \ref{intro}, and $R$ be the reflection of the direction in the vertical line; the latter takes $\alpha_{n_1\dots n_{k-1}n_k}$ to $\alpha_{3-n_1,\dots,3-n_{k-1},4-n_k}$. The action of $T_j$ on the symbolic orbits is described in Section \ref{intro}, while $R$ transforms the symbolic orbit $m_1\dots m_N$ to $6-m_1,\dots,6-m_N$. The transformation $(RT_4)^nRT_3$ takes $\alpha$ to $\alpha_{0\dots01}$ ($n-1$ zeroes); the transformation $T_j$ takes $\alpha_{0\dots01}$ to $\alpha_{j-1,3\dots33}$, and the transformation $R$ takes $\alpha_{j-1,3\dots33}$ to $\alpha_{4-j,0\dots01}$. This gives all the symbolic orbits listed above.
\proofend

{\bf  Induction step: short orbits}. 
We assume that, for some vertex of some pentagon, the short symbolic orbits satisfy our statement with some $x$ and $y$. Of these $x$ and $y$, we consider only the first and the last digit. We will see that for these digits only 6 possibilities occur:
$$
\begin{matrix}
(\1)&\quad&x=4\ast1&y=4\ast3\\
 (\2)&\quad&x=2\ast3&y=4\ast3\\
  (\3)&\quad&x=2\ast3&y=2\ast3\\
\end{matrix}
$$
and also $(\1)^\ast,(\2)^\ast,(\3)^\ast$ where $*$ means subtracting all the digits from 6.

First, we notice that, for $\alpha_1,\alpha_2,\alpha_3$, we have the cases $(\1),(\2),(\1)^\ast$. Then we will prove that the transformations $T_j,\, j=1,2,3,4,$ take the symbolic orbit $(xy)^nx, y(xy)^n$ with $x,y$ as in $({\bf k}),\, k=1,2,3,$ to symbolic orbit of the same form with $x$ and $y$ as in $({\bf l})$ or $({\bf l})^\ast$ (which we write as $T_j({\bf k})=({\bf l})$ or $({\bf l})^\ast$). More specifically,
$$
\begin{matrix}
T_1(\1)=(\3)\phantom{^\ast}&T_2(\1)=(\3)\phantom{^\ast}&T_3(\1)=(\3)^\ast&T_4(\1)=(\3)\phantom{^\ast}\\
 T_1(\2)=(\2)^\ast&T_2(\2)=(\3)\phantom{^\ast}&T_3(\2)=(\3)^\ast&T_4(\2)=(\3)\phantom{^\ast}\\
  T_1(\3)=(\3)^\ast&T_2(\3)=(\3)\phantom{^\ast}&T_3(\3)=(\3)^\ast&T_4(\3)=(\3)^\ast.\\
\end{matrix}
$$
Since  $T_j(({\bf k})^\ast)=(T_j({\bf k}))^\ast$, this will complete the induction.

In the computations below, we use, when necessary, Lemma \ref{lm}
\smallskip

\begin{tabular}{cll}
$(\1)$ &\qquad$(4\ast1\, 4\ast3)^n4\ast1$&\qquad \qquad \qquad \qquad$4\ast3\, (4\ast1\, 4\ast3)^n$\\
\end{tabular}
\smallskip
\hrule
\smallskip
\begin{tabular}{lll}
\hskip-.3in$T_1(\1)=(\3)$&\hskip-.12in$(\overrightarrow32\ast5(2)3\|2\ast2)^n3\|2\ast5(2)\langle3\rangle$&\hskip-.1in$\overrightarrow32\ast2(3\|2\ast5(2)3\|2\ast2)^n\langle3\rangle$\\
\hskip-.3in$T_2(\1)=(\3)$&\hskip-.12in$(2\ast4(3)\|2\ast1(43)\|)^n2\ast4(3)$&\hskip-.1in$2\ast1(43)\|(1\ast4(3)\|2\ast1(43))^n$\\
\hskip-.3in$T_3(\1)=(\3)^\ast$&\hskip-.12in$\langle4\rangle(1\ast3\|(4)1\ast5(23\|4))^n1\ast3(\overleftarrow4)$\\
\hskip-.3in&&\hskip-.4in$\langle4\rangle1\ast5(23\|4)(1\ast\|(4)1\ast5(23\overleftarrow4))^n$\\
\hskip-.3in$T_4(\1)=(\3)$&\hskip-.12in$\langle2\rangle(5\ast3\|25\ast4(3\|2))^n5\ast3\overleftarrow2$&\hskip-.3in$\langle2\rangle5\ast4(3\|2)(5\ast3\|25\ast4(3\|\overleftarrow2))^n$\\
\end{tabular}
\smallskip
\hrule
\smallskip
\begin{tabular}{lll}
$(\2)$&\qquad$(2\ast3\, 4\ast3)^n2\ast3$&\qquad \qquad \qquad \qquad \ \ $4\ast3\, (2\ast3\, 4\ast3)^n$\\
\end{tabular}
\smallskip
\hrule
\smallskip
\begin{tabular}{lll}
\hskip-.3in$T_1(\2)=(\2)^\ast$&\hskip-.1in$\langle4\rangle(1\ast23\|2\ast2(3\|4))^n1\ast2(3\overleftarrow4)$&\hskip-.3in$\overrightarrow32\ast2((3\|4)1\ast23\|2\ast2)^n\langle3\rangle$\\
\hskip-.3in$T_2(\2)=(\3)$&\hskip-.1in$\langle2\rangle(5\ast1(43)\|2\ast1(43\|2))^n5\ast1(43\overleftarrow2)$\\
\hskip-.3in&&\hskip-.45in$2\ast1(43\|2)(5\ast1(43)\|2\ast1(43))^n$\\
\hskip-.3in$T_3(\2)=(\3)^\ast$&\hskip-.1in$(4\ast5(23\|4)1\ast5(23)\|)^n4\ast5(23)$\\
&&\hskip-.8in$\langle4\rangle1\ast5((23)\|4\ast5(23\|4)1\ast5(23))^n\overleftarrow4$\\
\hskip-.3in$T_4(\2)=(\3)$&\hskip-.1in$(\overrightarrow32\ast4(3\|2)5\ast4)^n3\|2\ast4\langle3\rangle$&\hskip-.45in$\langle2\rangle5\ast4(3\|2\ast4(3\|2)5\ast4)^n(3\overleftarrow2)$\\
\end{tabular}
\smallskip
\hrule
\smallskip
\begin{tabular}{lll}
$(\3)$&\qquad$(\2\ast3\, 2\ast3)^n2\ast3$&\qquad \qquad \qquad \qquad \ \ $2\ast3\, (2\ast3\, 2\ast3)^n$\\
\end{tabular}
\smallskip
\hrule
\smallskip
\begin{tabular}{lll}
\hskip-.3in$T_1(\3)=(\3)^\ast$&\hskip-.1in$\langle4\rangle(1\ast2(3\|4)1\ast2(3\|4))^n1\ast2(3\overleftarrow4)$\\
&&\hskip-.7in$\langle4\rangle1\ast2(3\|4)(1\ast2(3\|4)1\ast2(3\overleftarrow4))^n$\\
\hskip-.3in$T_2(\3)=(\3)$&\hskip-.1in$\langle2\rangle(5\ast1(43\|2)5\ast1(43\|2))^n5\ast1(43\overleftarrow2)$\\
&&\hskip-.9in$\langle2\rangle5\ast1(43\|2)(5\ast1(43\|2)5\ast1)^n(43\overleftarrow2)$\\
\hskip-.3in$T_3(\3)=(\3)^\ast$&\hskip-.1in$(4\ast5(23)\|4\ast5(23)\|)^n4\ast5(23)$&\hskip-.45in$4\ast5(23)\|(4\ast5(23)\|4\ast5(23))^n$\\
\hskip-.3in$T_4(\3)=(\3)^\ast$&\hskip-.1in$\langle4\rangle(3\ast3\|43\ast3\|4)^n3\ast3\overleftarrow4$&\hskip-.15in$\langle4\rangle3\ast3\|4(3\ast3\|43\ast3\overleftarrow4)^n$\\
\end{tabular}
\smallskip

This completes the proof for short orbits.
\medskip

{\bf Induction step: long orbits}. 
For long symbolic orbits, we need to consider four fragments of sequences: $u,v,z,w$. According to the inductive assumption, the long orbits have the form $(zw)^nzv(uv)^n, (uv)^nuw(zw)^n$, and we apply the transformations $T_j$ to these orbits. It turns out that, in terms of  first and last digits, we will need to consider only 8 possibilities for $u,v,z,w$:
$$
\begin{matrix}
(\1)&u=2\ast3&v=\emptyset&z=2\ast3&w=4\ast3\\
 (\2)&u=4\ast3&v=2\ast3&z=4\ast3&w=2\ast3\\
  (\3)&u=2\ast3&v=\emptyset&z=2\ast3&w=2\ast3\\
   (\4)&u=2\ast3&v=2\ast3&z=2\ast3&w=2\ast3\\
\end{matrix}
$$
and also $(\1)^\ast,(\2)^\ast,(\3)^\ast,(\4)^\ast$.

Again, for $\alpha_1,\alpha_2,\alpha_3$, we have the cases $(\1),(\2),(\1)^\ast$. Then we consider the action of the transformations $T_j$ and prove the following result:
$$
\begin{matrix}T_1(\1)=(\1)^\ast&T_2(\1)=(\3)&T_3(\1)=(\3)^\ast&T_4(\1)=(\3)^\ast\\
 T_1(\2)=(\2)^\ast&T_2(\2)=(\4)&T_3(\2)=(\4)^\ast&T_4(\2)=(\2)^\ast\\
  T_1(\3)=(\3)^\ast&T_2(\3)=(\3)&T_3(\3)=(\3)^\ast&T_4(\3)=(\3)^\ast\\
   T_1(\4)=(\4)^\ast&T_2(\4)=(\4)&T_3(\4)=(\4)^\ast&T_4(\4)=(\4)^\ast\\
\end{matrix}
$$
This completes the induction for long orbits. The computations are shown below.
\smallskip

\begin{tabular}{lll}
$(\1)$&$(2\ast3\, 4\ast3)^n2\ast3\ \emptyset\ (2\ast3\ \emptyset\ )^n$&$(2\ast3\ \emptyset\ )^n2\ast3\, 4\ast3(2\ast3\, 4\ast3)^n$\\\end{tabular}
\smallskip
\hrule
\smallskip
\begin{tabular}{lll}
\hskip-.3in$T_1(\1)=(\1)^\ast$&\hskip-.1in$\langle4\rangle(1\ast23\|2\ast2(3\|4))^n1\ast2(3\|4)\ \emptyset\ (1\ast2(3\overleftarrow4)\ \emptyset\ )^n$\\
&&\hskip-2.95in$\langle4\rangle(1\ast2(3\|4)\ \emptyset\ )^n1\ast23\|2\ast2(3\|4)(1\ast23\|2\ast2(3\overleftarrow4))^n$\\
\hskip-.3in$T_2(\1)=(\3)$&\hskip-.1in$\langle2\rangle(5\ast1(43)\|2\ast1(43\|2))^n5\ast1(43\|2)\ \emptyset\ (5\ast1(43\overleftarrow2)\ \emptyset\ )^n$\\
&&\hskip-3.5in$\langle2\rangle(5\ast1(43)\|2\ \emptyset\ )^n5\ast1(43)\|2\ast1(43\|2)(5\ast1(43)\|2\ast1(43\overleftarrow2))^n$\\
\hskip-.3in$T_3(\1)=(\3)^\ast$&\hskip-.1in$(4\ast5(23\|4)1\ast5(23)\|)^n4\ast5(23)\|\ \emptyset\ (4\ast5(23)\ \emptyset\ )^n$\\
&&\hskip-3.2in$(4\ast5(23)\|\ \emptyset\ )^n4\ast5(23\|4)1\ast5(23)\|(4\ast5(23\|4)1\ast5(23))^n$\\
\hskip-.3in$T_4(\1)=(\3)^\ast$&\hskip-.1in$\langle4\rangle(3\ast3\|4(32)5\ast3\|4)^n3\ast3\|4\ \emptyset\ (3\ast3\overleftarrow4\ \emptyset\ )^n$\\
&&\hskip-2.9in$\langle4\rangle(3\ast3\|4\ \emptyset\ )^n3\ast3\|4(32)5\ast3\|4(3\ast3\|4(32)5\ast3\overleftarrow4)^n$\\
\end{tabular}
\smallskip
\hrule
\smallskip
\begin{tabular}{lll}
$(\2)$&\hskip1in$(4\ast3\, 2\ast3)^n4\ast3\, 2\ast3(4\ast3\, 2\ast3)^n$&\hskip.5in same\\
\end{tabular}
\smallskip
\hrule
\smallskip
\begin{tabular}{lll}
\hskip-.3in$T_1(\2)=(\2)^\ast$\\
&\hskip-.3in$(\overrightarrow32\ast2(3\|4)1\ast2)^n3\|2\ast2(3\|4)1\ast2(3\|2\ast2(3\|4)1\ast2)^n\langle3\rangle$&\hskip-.15in same\\
\hskip-.3in$T_2(\2)=(\4)$\\
&\hskip-.7in$(2\ast1(43\|2)5\ast1(43)\|)^n2\ast1(43\|2)5\ast1(43)\|(2\ast1(43\|2)5\ast1(43))^n$&\hskip-.15in same\\
\hskip-.3in$T_3(\2)=(\4)^\ast$\\
&\hskip-.85in$\langle4\rangle(1\ast5(23)\|4\ast5(23\|4))^n1\ast5(23)\|4\ast5(23\|4)(1\ast5(23)\|4\ast5(23\overleftarrow4))^n$&\hskip-.13in same\\
\hskip-.3in$T_4(\2)=(\2)^\ast$&\hskip-.1in$\langle2\rangle(5\ast43\|4\ast4(3\|3))^n5\ast43\|4\ast4(3\|2)(5\ast43\|4\ast4(3\overleftarrow2))^n$&\hskip-.13in same\\
\end{tabular}
\smallskip
\hrule
\smallskip
\begin{tabular}{lll}
$(\3)$&$\hskip.1in(2\ast3\, 2\ast3)^n2\ast3\ \emptyset\ (2\ast3\ \emptyset\ )^n$&$(2\ast3\ \emptyset\ )^n2\ast3\, 2\ast3(2\ast3\, 2\ast3)^n$\\\end{tabular}
\smallskip
\hrule
\smallskip
\begin{tabular}{lll}
\hskip-.3in$T_1(\3)=(\3)^\ast$&\hskip-.1in$\langle4\rangle(1\ast2(3\|4)1\ast2(3\|4))^n1\ast2(3\|4)\ \emptyset\ (1\ast2(3\overleftarrow4)\ \emptyset\ )^n$\\
&&\hskip-3.35in$\langle4\rangle(1\ast2(3\|4)\ \emptyset\ )^n1\ast2(3\|4)1\ast2(3\|4)(1\ast2(3\|4)1\ast2(3\overleftarrow4))^n$\\
\hskip-.3in$T_2(\3)=(\3)$&\hskip-.1in$\langle2\rangle(5\ast1(43\|2)5\ast1(43\|2))^n5\ast1(43\|2)\ \emptyset\ (5\ast1(43\overleftarrow2)\ \emptyset\ )^n$\\
&&\hskip-3.7in$\langle2\rangle(5\ast1(43\|2)\ \emptyset\ )^n5\ast1(43\|2)5\ast1(43\|2)(5\ast1(43\|2)5\ast1(43\overleftarrow2))^n$\\
\hskip-.3in$T_3(\3)=(\3)^\ast$&\hskip-.1in$(4\ast5(23)\|4\ast5(23)\|)^n4\ast5(23)\|\ \emptyset\ (4\ast5(23)\ \emptyset\ )^n$\\
&&\hskip-3.1in$(4\ast5(23)\|\ \emptyset\ )^n4\ast5(23)\|4\ast5(23)\|(4\ast5(23)\|4\ast5(23))^n$\\
\hskip-.3in$T_4(\3)=(\3)^\ast$&\hskip-.1in$(\overrightarrow34\ast43\|4\ast4)^n3\|4\ast4\ \emptyset\ (3\|4\ast4\ \emptyset\ )^n\langle3\rangle$\\
&&\hskip-2.5in$(\overrightarrow34\ast4\ \emptyset\ )^n3\|4\ast43\|4\ast4(3\|4\ast43\|4\ast4)^n\langle3\rangle$\\
\end{tabular}
\smallskip
\hrule
\smallskip
\begin{tabular}{lll}
$(\4)$&\hskip.3in$(2\ast3\, 2\ast3)^n2\ast3\, 2\ast3(2\ast3\, 2\ast3)^n$&\hskip1.23in same\\
\end{tabular}
\smallskip
\hrule
\smallskip
\begin{tabular}{lll}
\hskip-.3in$T_1(\4)=(\4)^\ast$\\
&\hskip-.7in$\langle4\rangle(1\ast3(3\|4)1\ast3(3\|4))^n1\ast3(3\|4)1\ast3(3\|4)(1\ast3(3\|4)1\ast3(3\overleftarrow4))^n$&\hskip-.5in same\\
\hskip-.3in$T_2(\4)=(\4)$\\
&\hskip-.7in$\langle2\rangle(5\ast1(43\|2)5\ast1(43\|2))^n5\ast1(43\|2)5\ast1(43\|2)(5\ast1(43\|2)5\ast1(43\overleftarrow2))^n$\\
&&\hskip-.5in same\\
\hskip-.3in$T_3(\4)=(\4)^\ast$\\
&\hskip-.7in$(4\ast5(23)\|4\ast5(23)\|)^n4\ast5(23)\|4\ast5(23)\|(4\ast5(23)\|4\ast5(23))^n$&\hskip-.5in same\\
\hskip-.3in$T_4(\4)=(\4)^\ast$&\hskip-.1in$\langle4\rangle(3\ast3\|43\ast3\|4)^n3\ast3\|43\ast3\|4(3\ast3\|43\ast3\overleftarrow4)^n$&\hskip-.5in same\\
\end{tabular}
\smallskip
\hrule
\medskip

This completes the proof for long orbits.

\bigskip

{\bf Acknowledgments}. D. F. is grateful to the MPIM Bonn for its hospitality. S. T. was partially supported by the Simons Foundation grant No 209361 and by the NSF grant DMS-1105442.

 \end{document}